\documentclass[12pt,reqno]{amsart}

\usepackage{latexsym}
\usepackage{amsfonts,amsmath,amssymb}
\usepackage{euscript}

\allowdisplaybreaks

\title[Elementary proof for a double sum]
{A short and elementary proof for a double sum of Brent and Osburn}

\author[H.~Prodinger]{Helmut Prodinger}
\address{Helmut Prodinger,
Mathematics Department, Stellenbosch University,
7602 Stellenbosch, South Africa.}
\email{hproding@sun.ac.za}

\begin{document}
\maketitle

We present here a completely elementary derivation of a recent formula of Brent and
Osburn~\cite{BrOs13}:
\begin{equation*}
S:=\sum_i\sum_j\binom{2n}{n+i}\binom{2n}{n+j}|i^2-j^2|=2n^2\binom{2n}{n}^2.
\end{equation*}
First, in order to get rid of the absolute value, we can rearrange the sum:
\begin{align*}
S&=2\sum_{i\ge0}\sum_{-i\le j\le i}\binom{2n}{n+i}\binom{2n}{n+j}(i^2-j^2)
+2\sum_{j\ge0}\sum_{-j\le i\le j}\binom{2n}{n+i}\binom{2n}{n+j}(j^2-i^2)\\
&=4\sum_{i\ge0}\sum_{-i\le j\le i}\binom{2n}{n+i}\binom{2n}{n+j}(i^2-j^2).
\end{align*}
Writing $i^2-j^2=-(n-i)(n+i)+(n-j)(n+j)$, we can continue:
\begin{align*}
\frac{S}{4(2n)(2n-1)}&=-\sum_{i\ge0}\sum_{-i\le j\le i}\binom{2n-2}{n-1+i}\binom{2n}{n+j}+
\sum_{i\ge0}\sum_{-i\le j\le i}\binom{2n}{n+i}\binom{2n-2}{n-1+j}\\
&=-2\sum_{0\le j\le i}\binom{2n-2}{n-1+i}\binom{2n}{n+j}
+\binom{2n}{n}\sum_{i\ge0}\binom{2n-2}{n-1+i}\\
&+2\sum_{0\le j\le i}\binom{2n}{n+i}\binom{2n-2}{n-1+j}-
\binom{2n-2}{n-1}\sum_{i\ge0}\binom{2n}{n+i}\\
&=-2\sum_{0\le j\le i}\binom{2n-2}{n-1+i}
\bigg[\binom{2n-2}{n+j}+2\binom{2n-2}{n-1+j}+\binom{2n-2}{n-2+j}\bigg]\\
&+2\sum_{0\le j\le i}\bigg[\binom{2n-2}{n+i}+2\binom{2n-2}{n-1+i}+\binom{2n-2}{n-2+i}\bigg]\binom{2n-2}{n-1+j}\\&-
\binom{2n-2}{n-1}\sum_{i\ge0}\binom{2n}{n+i}+\binom{2n}{n}\sum_{i\ge0}\binom{2n-2}{n-1+i}\\
&=-2\sum_{0\le j\le i}\binom{2n-2}{n-1+i}
\binom{2n-2}{n+j}+2\sum_{0\le j\le i}\binom{2n-2}{n-2+i}\binom{2n-2}{n-1+j}\\
&+2\sum_{0\le j\le i}\binom{2n-2}{n+i}\binom{2n-2}{n-1+j}-2\sum_{0\le j\le i}\binom{2n-2}{n-1+i}\binom{2n-2}{n-2+j}\\&-
\binom{2n-2}{n-1}\sum_{i\ge0}\binom{2n}{n+i}+\binom{2n}{n}\sum_{i\ge0}\binom{2n-2}{n-1+i}\\
&=2\binom{2n-2}{n-1}\sum_{ i\ge0}\binom{2n-2}{n-2+i}-2\binom{2n-2}{n-2}\sum_{i\ge0}\binom{2n-2}{n-1+i}\\&-
\binom{2n-2}{n-1}\sum_{i\ge0}\binom{2n}{n+i}+\binom{2n}{n}\sum_{i\ge0}\binom{2n-2}{n-1+i}\\
&=\frac{n}{4(2n-1)}\binom{2n}{n}^2.
\end{align*}
The last formula follows from the elementary 
\begin{equation*}
\sum_{i\ge0}\binom{2n}{n+i}=\frac12\bigg[4^n+\binom{2n}{n}\bigg].
\end{equation*}
It is not even necessary, since, with
\begin{equation*}
X:=\sum_{i\ge0}\binom{2n-2}{n-1+i},
\end{equation*}
the last expression equals
\begin{align*}
&2\binom{2n-2}{n-1}\bigg[\binom{2n-2}{n-2}+X\bigg]-2\binom{2n-2}{n-2}X\\&-
\binom{2n-2}{n-1}\bigg[4X-\binom{2n-2}{n-1}+\binom{2n-2}{n-2}\bigg]+\binom{2n}{n}X\\
&=\binom{2n-2}{n-1}\bigg[2\binom{2n-2}{n-2}+
\binom{2n-2}{n-1}-\binom{2n-2}{n-2}\bigg]=\binom{2n-2}{n-1}\binom{2n-1}{n-1}.
\end{align*}

\bibliographystyle{plain}

\end{document}